\documentclass[a4paper,10pt]{amsart}

\usepackage{amsmath,amssymb,amsthm,a4wide}

\usepackage{tikz}
\usetikzlibrary{shapes}
\usetikzlibrary{intersections}

\usepackage{graphicx}
\usepackage{graphics}

\newtheorem*{thm}{Theorem}

\newtheorem*{corollary}{Corollary}

\begin{document}

\title[QMC for Harmonic Functions]{A Remark on Disk Packings and Numerical\\ Integration of Harmonic Functions}
\author{Stefan Steinerberger}
\address{Stefan Steinerberger, Mathematisches Institut, Universit\"at Bonn, Endenicher Allee 60, 53115 Bonn, Germany}
\email{steinerb@math.uni-bonn.de}

\begin{abstract} We are interested in the following problem: given an open, bounded domain $\Omega \subset \mathbb{R}^2$, what is the largest constant $\alpha = \alpha(\Omega) > 0$ such that there exist an infinite sequence of disks $B_1, B_2, \dots, B_N, \dots \subset \mathbb{R}^2$ and a sequence $(n_i)$ with $n_i \in \left\{1,2\right\}$ such that
$$ \sup_{N \in \mathbb{N}}{N^{\alpha}\left\| \chi_{\Omega} - \sum_{i=1}^{N}{(-1)^{n_i}\chi_{B_i}}\right\|_{L^1(\mathbb{R}^2)}} < \infty,$$
where $\chi$ denotes the characteristic function?
We prove that certain (somewhat peculiar) domains $\Omega \subset \mathbb{R}^2$ satisfy the property with $\alpha = 0.53$. For these domains there exists a sequence of points $(x_i)_{i=1}^{\infty}$ in $\Omega$ with weights $(a_i)_{i=1}^{\infty}$ such that for all harmonic functions $u:\mathbb{R}^2 \rightarrow \mathbb{R}$
$$ \left|\int_{\Omega}{u(x)dx} - \sum_{i=1}^{N}{a_i u(x_i)}\right| \leq C_{\Omega}\frac{\|u\|_{L^{\infty}(\Omega)}}{N^{0.53}},$$
where $C_{\Omega}$ depends only on $\Omega$. This gives a Quasi-Monte-Carlo method for harmonic functions which improves on the probabilistic Monte-Carlo bound $\|u\|_{L^{2}(\Omega)}/N^{0.5}$ \textit{without} introducing a dependence on the total variation. We
do not know which decay rates are optimal.
\end{abstract}
\maketitle

\section{Introduction}
\subsection{Harmonic functions.} This paper aims to describe some progress in a problem that arose at the Oberwolfach Workshop 1340 'Uniform Distribution Theory and Applications', where it was motivated by a talk of the author on a related problem \cite{me}. We describe our question in its simplest possible setting: let $\Omega \subset \mathbb{R}^2$ be some bounded domain and let $u:\Omega \rightarrow \mathbb{R}$ be a harmonic function, i.e.
assume it satisfies 
$$ \Delta u = 0, \quad \mbox{where} \quad \Delta = \frac{\partial^2}{\partial x^2}+\frac{\partial^2}{\partial y^2}$$
is the Laplacian. Is there a Quasi Monte Carlo method able to exploit this information effectively to compute an approximation (including an error estimate) of
$$\int_{\Omega}{u(x)dx}?$$
The key ingredient suggesting that this might indeed be the case is the mean-value property: let $B(x,r)$ denote the disk with radius $r$ centered at $x \in \mathbb{R}^2$. If $u$ is harmonic in a neighbourhood of $B(x,r)$, then
$$ u(x) = \frac{1}{r^2\pi}\int_{B(x,r)}{u(z)dz}.$$
This means that exact integration over disks can be done with one function evaluation. In particular,
if one had a sequence of disks $B_i$ such that
$$ \sup_{N \in \mathbb{N}}{N^{\alpha}\left\| \chi_{\Omega} - \sum_{i=1}^{N}{(-1)^{n_i}\chi_{B_i}}\right\|_{L^1(\mathbb{R}^2)}}  \leq C_{\Omega}$$
for some $\alpha > 0$, then this gives a Quasi Monte Carlo method for harmonic functions
$$ \left|\int_{\Omega}{u(x)dx} - \sum_{i=1}^{N}{a_i u(x_i)}\right| \leq C_{\Omega}\frac{\|u\|_{L^{\infty}(\Omega)}}{N^{\alpha}},$$
where $x_i$ is the center of $B_i$ and $a_i = (-1)^{n_i}|B_i|$. Conversely, since the constant function 1 is harmonic any such Quasi Monte Carlo method gives a sequence of disks $B_i$ centered at $x_i$ with radius $r$ given via $r^2 \pi = a_i$ and $n_i = 1-(\mbox{sgn}(a_i)+1)/2$ such that 
$$ \sup_{N \in \mathbb{N}}{N^{\alpha}\left\| \chi_{\Omega} - \sum_{i=1}^{N}{(-1)^{n_i}\chi_{B_i}}\right\|_{L^1(\mathbb{R}^2)}} \leq C_{\Omega}.$$

\subsection{Main Result.} 
We will prove a result for the following (quite restricted but nontrivial) type of domains: we say $\Omega \subset \mathbb{R}^2$ is
\textit{finitely disk-covered} if there exists a finite number of closed disks $B_1,B_2, \dots, B_k$
such that any two disks meet at most in a single point that 'span' $\Omega$ in the following way: every point in $x \in \Omega$ is either contained in one of the disks or lies in a region surrounded by three disks such that any two out of these three disks touch in a point.

\begin{figure}[h!]
\begin{tikzpicture}[scale = 1 ]

\draw [ thick] (0,0) circle [radius=1];
\draw [ thick] (2,0) circle [radius=1];
\draw [ thick] (1,-0.65) circle [radius=0.2];
\draw [ thick] (1,0.8) circle [radius=0.29];

\draw [ultra thick,domain=35:325] plot ({-6+cos(\x)}, {sin(\x)});
\draw [ultra thick,domain=160:360] plot ({-5+0.2*cos(\x)}, {-0.65+0.2*sin(\x)});
\draw [ultra thick,domain=0:20] plot ({-5+0.2*cos(\x)}, {-0.65+0.2*sin(\x)});
\draw [ultra thick,domain=215:360] plot ({-4+cos(\x)}, {sin(\x)});
\draw [ultra thick,domain=0:145] plot ({-4+cos(\x)}, {sin(\x)});
\draw [ultra thick,domain=0:210] plot ({-5+0.29*cos(\x)}, {0.8+0.29*sin(\x)});
\draw [ultra thick,domain=330:360] plot ({-5+0.29*cos(\x)}, {0.8+0.29*sin(\x)});

\end{tikzpicture}
\caption{A simple example of a finitely disk-covered domain and the underlying disks.}
\end{figure}
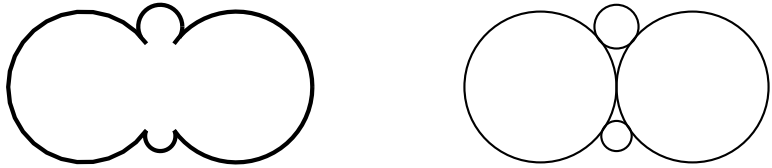
These sets are quite peculiar, however, at least any bounded, simply connected domain with a smooth boundary can be approximated in the Gromov-Hausdorff metric by a sequence of finitely disk-covered sets: it suffices
to consider disks in a lattice arrangement (either hexagonal or rectangular in which case on has to
add a final disk in the middle of every area) and approximate the desired domain using
this lattice (rescaled to the desired level of accuracy of the approximation).

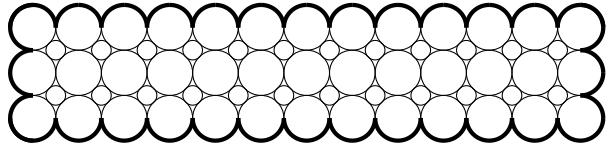
\begin{figure}[h!]
\begin{tikzpicture}[scale = 0.3]
 \foreach \j in {0,...,12}
       {\draw [ultra thick,domain=180:360] plot ({2*\j+cos(\x)}, {sin(\x)});} 

 \foreach \j in {0,...,1}
       {\draw [ultra thick,domain=270:360] plot ({24+cos(\x)}, {2+2*\j+sin(\x)});} 
 \foreach \j in {0,...,2}
       {\draw [ultra thick,domain=0:90] plot ({24+cos(\x)}, {0+2*\j+sin(\x)});} 

 \foreach \j in {0,...,12}
       {\draw [ultra thick,domain=0:90] plot ({2*\j+cos(\x)}, {4+sin(\x)});} 
 \foreach \j in {0,...,12}
       {\draw [ultra thick,domain=90:180] plot ({2*\j+cos(\x)}, {4+sin(\x)});} 

\foreach \j in {0,...,2}
       {\draw [ultra thick,domain=90:270] plot ({+cos(\x)}, {2*\j+sin(\x)});} 

 \foreach \j in {0,...,12}
       {\draw (2*\j,0) circle [radius=01];} 
 \foreach \j in {0,...,12}
       {\draw (2*\j,2) circle [radius=01];} 
 \foreach \j in {0,...,12}
       {\draw (2*\j,4) circle [radius=01];} 

 \foreach \j in {0,...,11}
       {\draw (2*\j+1,1) circle [radius=0.43];} 

 \foreach \j in {0,...,11}
       {\draw (2*\j+1,3) circle [radius=0.43];} 

\end{tikzpicture}
\caption{A simple approximation of a rectangle: circles of equal radius in a rectangular grid with smaller circles filling up the holes.}
\end{figure}

Our consideration of this particular class of sets is twofold: it effectively cuts off the possibility
of harmonic functions with large growth at the boundary of the domain since numerical integration
close to the boundary of a finitely disk-covered domain can be done with finitely many function
evaluations; secondly, the arising structure allows us to exploit recent advances in the study of
Apollonian packings.
\begin{thm} Let $\Omega$ be finitely disk-covered. Then there exists
a sequence of disks $B_i$ such that
$$ \sup_{N \in \mathbb{N}}{N^{0.53}\left\| \chi_{\Omega} - \sum_{i=1}^{N}{\chi_{B_i}}\right\|_{L^1(\mathbb{R}^2)}} \leq C_{\Omega}.$$
where $C_{\Omega}$ depends only on $\Omega$.
\end{thm}
 We believe that the statement is not optimal and that one should be able to construct sequences with a larger exponent in $N$ if one were to exploit the fact that some of the disks may have negative coefficients (something that is not used in the statement here). It should certainly be possible to prove some bounds for, say, the class of convex domains (numerical experiments suggest that the randomized greedy algorithm -- picking a random point and, if it is not contained in one the existing disks, add the largest disk possible without introducing intersections -- corresponds to $\alpha \sim 0.2$ for convex domains). Since we never actually use the
possibility of subtracting characteristic functions (i.e. $n_i = 1$), we can immediately deduce that
for any $p \geq 1$
$$ \sup_{N \in \mathbb{N}}{N^{\frac{0.53}{p}}\left\| \chi_{\Omega} - \sum_{i=1}^{N}{\chi_{B_i}}\right\|_{L^p(\mathbb{R}^2)}} \leq C_{\Omega}$$
but there is no reason to assume that this should be in any way optimal. 

\subsection{Quasi-Monte Carlo.} As outlined above, the statement immediately implies
a Quasi Monte Carlo method. The same considerations as above suggest that there is
no reason to assume this might be optimal.
 \begin{corollary} Let $\Omega$ be finitely disk-covered. Then there exists
a universal sequence $(x_i)_{i=1}^{\infty}$ of points in $\Omega$ and a sequence $(a_i)_{i=1}^{\infty}$ of nonnegative reals with the following property: if
$$ \Delta u = 0 \qquad \mbox{in a neighbourhood of}~\Omega,$$
then
$$ \left| \int_{\Omega}{u(x)dx} - \sum_{i=1}^{N}{a_i u(x_i)}\right| \leq C_{\Omega} \frac{ \|u\|_{L^{\infty}(\Omega)} }{N^{0.53}},$$
where $C_{\Omega}$ depends only on $\Omega$.
\end{corollary}
Let us emphasize the difference to classical QMC methods: a Quasi-Monte-Carlo method is based on the simple approximation
$$  \int_{[0,1]^2}{u(x)dx} \sim \frac{1}{N}\sum_{i=1}^{N}{u(x_i)}$$
for a set of points $(x_i)_{i=1}^{N}$. The well-known Koksma-Hlawka inequality gives
$$ \left| \int_{[0,1]^2}{u(x)dx} - \frac{1}{N}\sum_{i=1}^{N}{u(x_i)}\right| \leq D_N(x_i) V(f),$$
where $D_N$ denotes the discrepancy of the point set and $V(f)$ the total variation in the sense of Hardy-Krause. We refer to the classical monographs of Kuipers \& Niederreiter \cite{kui}, Drmota \& Tichy \cite{drm} and Dick \& Pillichshammer \cite{dick} for further information. We emphasize that there exist point sets such that $D_N(x_i) \sim N^{-1}\log{N}.$ However, and this is crucial, our bound is independent of the total variation of the function. Indeed, for the harmonic function (given
in polar coordinates)
$$ u_{m}(r, \theta) = r^m \cos{(m\theta)},$$
on some domain, we easily see that $V(u_m) \sim m$, which can be made arbitrarily large; in contrast, our bound is independent of $m$. 

\subsection{Possible extensions.} If we were to modify the approximation scheme using a suitably rescaling, then for suitable points and weights the approximation 
$$ \int_{\Omega}{u(x)dx} \sim \left(\frac{|\Omega|}{\sum_{i=1}^{N}{a_i}}\right)\sum_{i=1}^{N}{a_i u(x_i)}$$
should yield even better results: what decay properties can be proven? Another natural conjecture is that, at least for finitely disk-covered domains, even
$$ \left| \int_{\Omega}{u(x)dx} - \sum_{i=1}^{N}{a_i u(x_i)}\right| \leq C_{\Omega} \frac{ \|u\|_{L^{1}(\Omega)} }{N^{0.53}}$$
might be true. 

\section{The Proof}
\begin{proof}[Proof of the Theorem.]
The proof is constructive: since $\Omega$ is finitely disk-covered, we are initially given a 
finite set of disks $D_1, D_2, D_3, \dots, D_k$ associated to $\Omega$ with centers $x_1, \dots, x_k$. The mean-value theorem implies that for any harmonic $u$
$$ \sum_{i=1}^{k}{|D_i| u(x_i)} = \int_{\bigcup_{i=1}^{k}{D_i}}{u(x)dx}.$$
This is already precise on some part of the domain. The idea is to cover the rest of the domain with smaller and smaller disks (on each of which exact integration can again be performed). Let us consider a connected component of
$$ \Omega \setminus \bigcup_{i=1}^{k}{D_i}.$$
By assumption, it is bounded by three disks any two of which mutually touch in a point. Then
there exists precisely one circle contained within the connected domain that is tangent to all
three boundary circles: the statement dates back to Apollonius. Such a configuration of 4 circles
 is known as a Descartes configuration: given
a Descartes configuration, it is possible to construct three additional circles within the three
gaps. Iterating this process yields an Apollonian packing.
\begin{figure}[h!]
\begin{tikzpicture}[scale = 5]
\draw [ultra thick,domain=110:180] plot ({-1+1.8+0.8*cos(\x)}, {0.8*sin(\x)});
\draw [ultra thick,domain=0:65] plot ({-1+cos(\x)}, {sin(\x)});
\draw [ultra thick,domain=230:295] plot ({-1+1.113+cos(\x)}, {1.672+sin(\x)});
\draw [ultra thick,domain=0:360] plot ({-1+1.017+0.144*cos(\x)}, {0.535+0.144*sin(\x)});

\draw [ultra thick,domain=110:180] plot ({1.8+0.8*cos(\x)}, {0.8*sin(\x)});
\draw [ultra thick,domain=0:65] plot ({cos(\x)}, {sin(\x)});

\draw [ultra thick,domain=230:295] plot ({1.113+cos(\x)}, {1.672+sin(\x)});
\draw [ultra thick,domain=0:360] plot ({1.017+0.144*cos(\x)}, {0.535+0.144*sin(\x)});
\draw [ultra thick,domain=0:360] plot ({0.84+0.06*cos(\x)}, {0.643+0.06*sin(\x)});
\draw [ultra thick,domain=0:360] plot ({1.2+0.06*cos(\x)}, {0.62+0.06*sin(\x)});
\draw [ultra thick,domain=0:360] plot ({1.005+0.058*cos(\x)}, {0.335+0.058*sin(\x)});
\end{tikzpicture}
\caption{Left: a Descartes configuration. Right: adding three additional circles}
\end{figure}
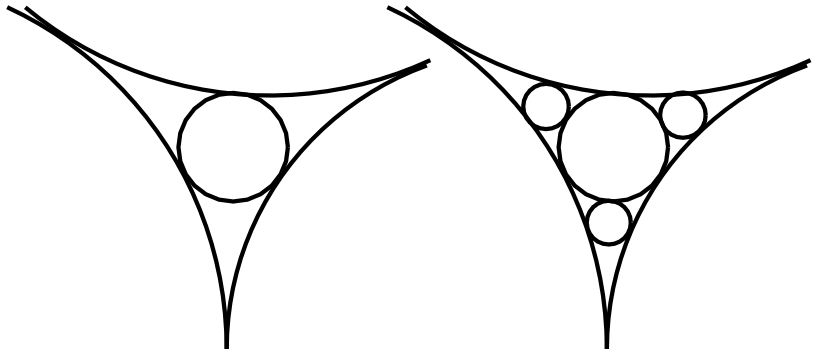
For each connected component of $ \Omega \setminus \bigcup_{i=1}^{k}{D_i}$ (of which there are only finitely many) we construct the associated Apollonian packing and then define an infinite sequence of disks $E_1, E_2, \dots$
by ordering the union of the disks created by the Apollonian packings and the finitely many disks $D_1, D_2,\dots, D_k$ by size. Let $x_i$ denote the center of $E_i$. Using the mean-value property, we get that
$$ \left|  \sum_{i=1}^{N}{|E_i| u(x_i)} - \int_{\Omega}{u(x)dx}\right| \leq
 \left|\int_{\Omega \setminus \bigcup_{i=1}^{N}{E_i}}{u(x)dx}\right| \leq \left|\Omega \setminus \bigcup_{i=1}^{N}{E_i}\right| \|u\|_{L^{\infty}(\Omega)}.$$
It remains to control the speed with which the disks exhaust the set. Here we use a recent result
of Kontorovich \& Oh \cite{kant}: generalizing an earlier result of Boyd \cite{boy}, they show
the cardinality of disks with curvature $\kappa$ bounded from above by $T$ behaves as
$$ c_1 \cdot T^{\alpha} \leq \#\left\{i \in \mathbb{N}: \kappa(E_i) \leq T\right\} \leq c_2\cdot T^{\alpha}$$
for a universal constant $\alpha \sim 1.30568\dots$ (this approximation is due to McMullen \cite{mc}) and constants $c_1, c_2$ depending on the particular Apollonian packing. We consider merely finite number of Apollonian packings at the same time and may thus fix the constants $c_1, c_2$ in what follows. This implies that
$$c_2T^{\alpha} \leq \#\left\{i \in \mathbb{N}: T \leq \kappa(E_i) \leq \left(\frac{2c_2}{c_1}
\right)^{\frac{1}{\alpha}}T\right\} \leq \left(\frac{2c_2^2}{c_1}-c_1\right)T^{\alpha} \qquad \quad (\diamond).$$
This estimate controls the number of disks with curvature in a certain interval and shows that for on average there are $\sim T^{\alpha - 1}$ disks with curvature $T \leq \kappa \leq T+1$. We have that
$$  \left|\Omega \setminus \bigcup_{i=1}^{\infty}{E_i}\right| = 0$$
and therefore
$$ \left|\Omega \setminus \bigcup_{\kappa(E_i) \leq T}^{\infty}{E_i}\right| = 
 \left|\bigcup_{\kappa(E_i) \geq T}^{\infty}{E_i}\right|.$$
A disk with curvature $\kappa$ has measure $\pi/\kappa^2$ and therefore using $(\diamond)$
\begin{align*} \left|\Omega \setminus \bigcup_{\kappa(E_i) \leq T}^{}{E_i}\right| &\leq
\sum_{n=0}^{\infty}{\frac{\pi}{\left(\frac{2c_2}{c_1}
\right)^{\frac{2n}{\alpha}}T^2} \#\left\{i \in \mathbb{N}:  \left(\frac{2c_2}{c_1}
\right)^{\frac{n}{\alpha}}T \leq \kappa(E_i) \leq \left(\frac{2c_2}{c_1}
\right)^{\frac{n+1}{\alpha}}T\right\}} \\
&\leq \sum_{n=0}^{\infty}{\frac{\pi}{\left(\frac{2c_2}{c_1}
\right)^{\frac{2n}{\alpha}}T^2} \left(\frac{2c_2^2}{c_1}-c_1\right)\left(\frac{2c_2}{c_1}
\right)^{n}T^{\alpha}}\\
&\leq c\cdot T^{\alpha -2},
\end{align*}
for some constant $c$. If we define $N$ to be the number of circes with curvature bounded from above by $T$, then
$$ N \sim T^{\alpha} \quad \mbox{and thus} \qquad T^{\alpha - 2} \sim
N^{\frac{\alpha-2}{\alpha}} \sim  \frac{1}{N^{\frac{2-\alpha}{\alpha}}}.$$
Since $\alpha \sim 1.30568\dots$, we have that
$$ \frac{2-\alpha}{\alpha} = 0.536\dots$$
and this yields the result. \end{proof}

\section{Open problems} 
\subsection{Optimal decay rates.} The natural question is which decay rates are optimal. Our proof may be regarded as a greedy algorithm: the big open
question is the following: given a domain $\Omega \subset \mathbb{R}^2$, is it true
that the best approximation of $\chi_{\Omega}$ is always given by
$$ \chi_{\Omega} \sim \chi_{B_1} + \chi_{B_2} + \dots + \chi_{B_n}$$
for a sequence of balls $(B_i)$ or whether there exist more interesting configurations for which
$$ \chi_{\Omega} \sim \pm \chi_{B_1} \pm \chi_{B_2} \pm \dots \pm \chi_{B_n}$$
yields a better result for a suitable choice of signs. 

\subsection{Harmonic functions on fractal sets.} We conjecture that on finitely disk-covered domains for the sequence of disks constructed in the argument and an arbitrary harmonic function $u$ actually the following stronger inequality should be true
$$ \left| \int_{\Omega}{u(x)dx} - \sum_{i=1}^{N}{a_i u(x_i)}\right| \leq C_{\Omega} \frac{ \|u\|_{L^{1}(\Omega)} }{N^{0.53}}.$$
We emphasize that this is not a geometric statement about the constructed packing of disks and that
the statement is trivially false for arbitrary functions $u$. Our reasoning behind conjecturing such an
inequality is the fact that the set
$$ \Omega \setminus \bigcup_{i_1}^{N}{B_i} \qquad \mbox{has a very fractal structure}$$
while harmonic functions have strong ridigity properties. It seems extremely natural to assume that
harmonic functions cannot differ too much on fractal sets from their average behavior.\\

\textbf{Acknowledgement.} I am grateful to Michael Gnewuch and the other organizers of the Oberwolfach Workshop 1340 for the opportunity to present my work and all participants for the very enjoyable week. I am also happy to acknowledge an interesting discussion with Stefan Heinrich after having presented the problem. This work was supported by SFB 1060 of the DFG.

\end{document}